\documentclass[a4paper,12pt]{article}

\usepackage{latexsym}
\usepackage{mathrsfs}
\usepackage{amsfonts,amsmath,amssymb}
\usepackage{indentfirst}
\usepackage{subeqnarray}
\usepackage{verbatim}
\usepackage[pdftex]{graphicx}

\usepackage{epstopdf}
\usepackage{booktabs}
\usepackage{geometry}
\usepackage{bm}
\usepackage{enumerate}
\usepackage{float}

\usepackage{algorithmicx}
\usepackage[ruled]{algorithm}
\usepackage{algpseudocode}

\newtheorem{theorem}{Theorem}
\newtheorem{lemma}[theorem]{Lemma}

\numberwithin{equation}{section}
\newenvironment{pfs}[1][Proof]{\noindent\textbf{#1.} }
{\ \rule{0.75em}{0.75em}\smallskip}

\DeclareMathOperator*{\argmin}{arg\,min}

\begin{document}


\begin{center}
\Large\bf A nonsmooth optimization approach\\ for hemivariational inequalities \\ with applications in Contact Mechanics
\end{center}

\begin{center}
Michal Jureczka\footnote{Jagiellonian University in Krakow, Faculty of Mathematics and Computer Science, 
Lojasiewicza 6, 30-348 Krakow, Poland. Email: {\tt michal.jureczka@uj.edu.pl}}\quad and \quad 
Anna Ochal\footnote{Jagiellonian University in Krakow, Faculty of Mathematics and Computer Science, 
Lojasiewicza 6, 30-348 Krakow, Poland. Email: {\tt anna.ochal@uj.edu.pl}}
\end{center}

\begin{quote}
{\bf Abstract}. In this paper we introduce an abstract nonsmooth optimization problem and prove
existence and uniqueness of its solution. We present a numerical scheme to approximate this solution. The theory is later applied to a sample static contact problem describing an elastic body in frictional contact with a foundation. This contact is governed by a nonmonotone friction law with dependence on normal and tangential components of displacement. Finally, computational simulations are performed to illustrate obtained results.

{\bf Keywords}. Nonmonotone friction, optimization problem, error estimate, finite element method, numerical simulations.

{\bf AMS Classification.} 35Q74, 49J40, 65K10, 65M60, 74S05, 74M15, 74M10, 74G15
\end{quote}


\section{Introduction}
\noindent
In the literature we can find examples of many models describing displacement of deformable body that is partly in contact with another object, the so-called foundation. In various contact models boundary conditions enforced on the part of the body contacting the foundation appear. Functions that occur in these conditions model response of the foundation in direction normal to the contact boundary and in direction tangential to the boundary (friction law). In many cases these functions are monotone, such as when Coulomb's law of dry friction is considered, but in applications this may not always be the case. What is more, the friction bound may change as the penetration of the foundation by body increases. Nonmonotonicity of functions describing contact laws and influence of normal displacement of the body on friction law cause some difficulties in analytical and numerical treatment of considered problems.

In this paper we introduce an abstract framework that can be used to numerically approximate a solution to a class of mechanical contact problems. We present a nonsmooth optimization problem and prove existence and uniqueness of a solution to this problem. Next we present~a numerical scheme approximating this solution and provide numerical error estimation. We apply this theory to a static contact problem describing an elastic body in contact with a foundation. This contact is governed by a nonmonotone friction law with dependence on normal and tangential components of displacement. Weak formulation of introduced contact problem is presented in the form of hemivariational inequality. In the end we show results of computational simulations and describe the numerical algorithm that was used to obtain these results.

Let us now briefly present references in the literature. The definition and properties of Clarke subdifferential and tools used to solve optimization problems were introduced in $\cite{C}$. Comparison of nonsmooth and nonconvex optimization methods can be found in $\cite{BKM}$, and details on computational contact mechanics is presented in $\cite{W}$. The theory of hemivariational inequalities was developed in $\cite{P}$, and the idea to use Finite Element Method to solve these inequalities was presented in $\cite{HMP}$. Another early study of vector-valued hemivariational problems in the context of FEM can be found in $\cite{MH}$. More recent analysis of hemivariational and variational-hemivariational inequalities was presented in $\cite{MOS2}$, $\cite{MOS}$, whereas numerical analysis of such problems can be found for example in papers $\cite{BBK}$, $\cite{BBKR}$, $\cite{BHM}$, $\cite{H}$, $\cite{HSB}$, $\cite{HSD}$.

A similar mechanical model to the one described in the paper was already considered in $\cite{MOS}$, where the authors prove only existence of a solution using surjectivity result for pseudomonotone, coercive multifunction without requiring any smallness assumption.

An error estimation concerning stationary variational-hemivariational inequalities was presented in $\cite{H}$. In our case variational part of inequality is not present and the inequality is not constrained, however error estimations had to be generalized to reflect dependence of friction law on normal component of the displacement. 

A numerical treatment of mechanical problem leading to hemivariational inequality using two approaches - nonsmooth and nonconvex optimization and quasi-augmented Lagrangian method is presented in $\cite{BBK}$. As the smallness assumption is not required, this once again does not guarantee uniqueness and leads to a nonconvex optimization problem. There, the authors assume contact to be bilateral and consider friction law which does not depend on normal component of the displacement.

This paper is organized as follows. Section~$\ref{op}$ contains a general differential inclusion pro\-blem and an optimization pro\-blem. We show that under introduced assumptions both problems are equivalent and have a unique solution. In Section~$\ref{ns}$ we proceed with a discrete scheme that approximates solution to introduced optimization problem and we prove theorem concerning numerical error estimation. An application of presented theory in the form of mechanical contact model is indicated in Section~$\ref{mcp}$, along with its weak formulation. Finally, in Section~$\ref{s}$, we describe computational algorithm used to solve mechanical contact problem and present simulations for a set of sample data.

\section{A general optimization problem} \label{op}
\noindent
Let us start with basic notation used in this paper. For a normed space $X$, we denote by $\|\cdot\|_X$ its norm, by $X^*$ its dual space and by $\langle\cdot,\cdot\rangle_{X^*\times X}$ the duality pairing of $X^*$ and $X$. By $c>0$ we denote a generic constant (value of $c$ may
differ in different equations).

Let us now assume that $j\colon X \to \mathbb{R}$ is locally Lipschitz continuous. The generalized directional derivative of $j$ at $x \in X$ in the direction $v\in X$ is defined by
\begin{align*}
  &j^0(x;v) := \limsup_{y \to x, \lambda \searrow 0} \frac{j(y+\lambda v)  - j(y)}{\lambda}.
\end{align*}
The generalized subdifferential of $j$ at $x$ is a subset of the dual space $X^*$ given by
\begin{align*}
  &\partial j(x) := \{\xi \in X^*\, | \, \langle \xi, v\rangle_{X^*\times X}\leq j^0(x;v) \  \mbox{ for all } v \in X \}.
\end{align*}
If $j\colon X^n \to \mathbb{R}$ is a locally Lipschitz function of $n$ variables, then we denote by $\partial_i j$ and $j_i^0$ the Clarke subdifferential and generalized directional derivative with respect to $i$-th variable of~$j$, respectively.

Let now $V$ be a reflexive Banach space and $X$ be a Banach space. Let $\gamma \in \mathcal{L}(V, X)$ be linear and continuous operator from $V$ to $X$, and $c_{\gamma}:=\|\gamma\|_{\mathcal{L}(V, X)}$. We denote by  $\gamma^* \colon X^* \to V^*$ the adjoint operator to $\gamma$.  Let $A \colon V \to V^*$, $J \colon X \times X \to \mathbb{R}$ and $f\in V^*$. We formulate the differential inclusion problem as follows.

\medskip
\noindent
\textbf{Problem $\bm{P_{incl}}$:} {\it\ Find $u \in V$ such that}
\begin{align*}
  &A u + \gamma^* \partial_2 J(\gamma u, \gamma u) \ni f. 
\end{align*}

In the study of Problem $P_{incl}$ we make the following assumptions.

\medskip
\noindent
$\underline{H(A)}:$ \quad The operator $A \colon V \to V^*$ is such that
\begin{enumerate}
  \item[(a)]
     $A$ is linear and bounded,
  \item[(b)]
    $A$ is symmetric, i.e. $\langle A u, v\rangle_{V^*\times V}  = \langle A v ,  u  \rangle_{V^*\times V}$ for all $u, v \in V$,
  \item[(c)]
     there exists $m_A>0$ such that $\langle A u, u\rangle_{V^*\times V} \geq m_A \|u\|_V^2$ for all $u \in V$.
\end{enumerate}

\medskip
\noindent
$\underline{H(J)}:$ \quad The functional $J \colon X \times X \to \mathbb{R}$ satisfies
\begin{enumerate}
  \item[(a)]
     $J$ is locally Lipschitz continuous with respect to its second variable,
  \item[(b)]
     there exist $c_{0}, c_{1}, c_{2} \geq 0$ such that \\[2mm]
     \hspace*{1cm}$\|\partial_2 J(w, v)\|_{X^*} \leq c_{0} + c_{1}\| v \|_X + c_{2}\| w \|_X$\quad for all $w, v \in X$,
  \item[(c)]
    there exist $m_\alpha, m_L  \geq 0$ such that \\[2mm]
    $J_2^0(w_1, v_1; v_2 - v_1) + J_2^0(w _2, v_2; v_1 - v_2)\leq m_\alpha\| v_1- v_2\|_X^2+ m_L\| w_1 -  w_2 \|_{X} \| v_1 -  v_2\|_{X}$ \\[2mm]
    for all $ w_1,  w_2,  v_1,  v_2 \in X$.
\end{enumerate}

\medskip
\noindent
$\underline{H(f)}: \quad f \in V^*$.

\medskip
\noindent
$\underline{(H_s)}: \quad m_A > (m_\alpha + m_L) c_\gamma^2$.

\bigskip
We remark that condition $H(J)$(c) is a more general form of a relaxed monotonicity condition, i.e. for all  $w_1, w_2, v_1, v_2 \in X$
\begin{align*}
 &\langle \partial_2 J(w_1,  v_1) - \partial_2 J(w_2,  v_2), v_1 -  v_2 \rangle_{X^*\times X} \geq - m_\alpha\| v_1- v_2\|_X^2 - m_L\| w_1- w_2\|_X \| v_1- v_2\|_X.
\end{align*}
Moreover, in a special case when $J$ does not depend on the first variable (i.e. $w_1=w_2=w$), condition $H(J)$(c) is equivalent to a relaxed monotonicity condition, i.e. for all  $w, v_1, v_2 \in X$
\begin{align}
 &\langle \partial_2 J(w,  v_1) - \partial_2 J(w,  v_2), v_1 -  v_2 \rangle_{X^*\times X} \geq - m_\alpha \| v_1- v_2\|_X^2.  \label{RM}
\end{align}

\noindent
We start with a uniqueness result for Problem~$P_{incl}$.
\begin{lemma} \label{I}
Assume that $H(A)$, $H(J)$, $H(f)$ and $(H_s)$ hold. If Problem~$P_{incl}$ has a~solution $u \in V$, then it is unique and satisfies
\begin{align}
  &\| u \|_V \leq c \, (1+\| f \|_{V^*}) \label{e2}
\end{align}
with a positive constant $c$.
\end{lemma}

\begin{pfs}
Let $ u  \in V$ be a solution to Problem~$P_{incl}$. This means that there exists $ z  \in \partial_2 J(\gamma u , \gamma u )$ such that
\begin{align*}
  &A u  + \gamma^*  z  =  f .
\end{align*}
From the definition of generalized directional derivative of $J(\gamma u, \cdot)$ we have for all $ v  \in V$
\begin{align}
  &\langle  f  -  A u ,  v  \rangle_{V^*\times V} = \langle \gamma^* z ,  v \rangle_{V^*\times V} = \langle  z , \gamma v \rangle_{X^*\times X} \leq J_2^0(\gamma u , \gamma u ;  \gamma v ). \label{9}
\end{align}
Let us now assume that Problem $P_{incl}$ has two different solutions $u_1$ and $u_2$. For a~solution $u_1$ we set $ v  =  u _2 -  u _1$ in ($\ref{9}$) to get
\begin{align*}
  &\langle  f ,  u _2- u _1 \rangle_{V^*\times V} - \langle A u _1,  u _2 -  u _1 \rangle_{V^*\times V} \leq J_2^0(\gamma u _1, \gamma u _1; \gamma u _2-\gamma u _1).
\end{align*}
For a solution $u_2$ we set  $ v  =  u _1 -  u _2$ in ($\ref{9}$) to get
\begin{align*}
  &\langle  f ,  u _1- u _2 \rangle_{V^*\times V} - \langle A u _2,  u _1 -  u _2 \rangle_{V^*\times V} \leq J_2^0(\gamma u _2, \gamma u _2; \gamma u _1-\gamma u _2).
\end{align*}
Adding the above inequalities, we obtain
\begin{align*}
  &\langle A u _1 - A u _2,  u _1 -  u _2 \rangle_{V^*\times V}\\
  &\quad \leq  J_2^0(\gamma u _1, \gamma u _1; \gamma u _2-\gamma u _1) + J_2^0(\gamma u _2, \gamma u _2; \gamma u _1-\gamma u _2).
\end{align*}
Hence, $H(A)$(c) and $H(J)$(c) yield
\begin{align*}
  &m_A \|  u _1 -  u _2\|_V^2 \leq (m_\alpha + m_L)\|\gamma  u _1-\gamma  u _2\|_X^2,
\end{align*}
and finally
\begin{align*}
  &\big(m_A - (m_\alpha + m_L)c_\gamma^2\big) \|  u _1 -  u _2\|_V^2 \leq 0.
\end{align*}
Under assumption $(H_s)$, we obtain that if Problem $P_{incl}$ has a solution, it is unique.

Now, in order to prove $(\ref{e2})$, we set $ v  = - u $ in ($\ref{9}$) to obtain
\begin{align}
  &\langle A u ,  u  \rangle_{V^*\times V} \leq J_2^0(\gamma  u , \gamma  u ; -\gamma  u ) + \langle  f ,  u  \rangle_{V^*\times V}. \label{111}
\end{align}
Using $H(J)$(b) and (c), we get
\begin{align}\label{112}
&J_2^0(\gamma  u , \gamma  u ;-\gamma  u ) \leq (m_\alpha + m_L) \|\gamma  u \|_X^2 - J_2^0(0 , 0 ; \gamma  u ) \nonumber\\
&\quad \leq  (m_\alpha + m_L) \|\gamma  u \|_X^2 + c_{0} \|\gamma  u \|_X.
\end{align}
Combining ($\ref{111}$) and ($\ref{112}$), we have
\begin{align*}
m_A\| u \|_V^2 &\leq (m_\alpha + m_L) \|\gamma  u \|_X^2 + c_{0}\| \gamma  u \|_X + \| f \|_{V^*}\| u \|_V
\end{align*}
and
\begin{align*}
\big(m_A - (m_\alpha + m_L) c_\gamma^2\big)\| u \|_V \leq  c \, (1+\| f \|_{V^*}).
\end{align*}
From $(H_s)$ we obtain required estimation.
\end{pfs}

We now consider an optimization problem, which will be equivalent to Problem~$P_{incl}$ under introduced assumptions.
To this end, let the operator $\mathcal{L}: V \times V \rightarrow \mathbb{R}$ be defined for all $ w ,  v  \in V$ as follows
\begin{align}
  \mathcal{L}( w ,  v ) = \frac{1}{2} \langle A  v ,  v  \rangle_{V^*\times V}  - \langle  f ,  v  \rangle_{V^*\times V} + J(\gamma  w  ,\gamma  v  ). \label{L}
\end{align}
The next lemma collects some properties of the operator~$\mathcal{L}$.
\begin{lemma} \label{Lprop}
Under assumptions  $H(A)$, $H(J)$, $H(f)$ and $(H_s)$, the operator\\ $\mathcal{L}\colon V \times V \rightarrow \mathbb{R}$ defined by {\rm{(\ref{L})}} satisfies
\begin{enumerate}[(i)]
 \item $\mathcal{L}(w,\cdot)$ is locally Lipschitz continuous for all $w \in V$,
 \item $\partial_2 \mathcal{L}( w ,  v ) \subseteq A  v  -  f  + \gamma^* \partial_2 J(\gamma  w , \gamma  v )$ for all $w,v \in V$,
 \item $\mathcal{L}(w,\cdot)$ is strictly convex for all $w \in V$.
\end{enumerate}

\end{lemma}

\begin{pfs}
The proof of $(i)$ is immediate since for a fixed $w \in V$ the operator $\mathcal{L}(w,\cdot)$ is locally Lipschitz continuous as a sum of locally Lipschitz continuous functions with respect to $v$.

\medskip
\noindent
For the proof of $(ii)$, we observe that from $H(A)$ and $H(f)$, the functions
\begin{align*}
f_1\colon V \ni  v \mapsto \frac{1}{2} \langle A v ,  v  \rangle_{V^*\times V} \in \mathbb{R}, \qquad  f_2\colon V \ni v \mapsto \langle  f ,  v  \rangle_{V^*\times V} \in \mathbb{R}
\end{align*}
are strictly differentiable and we calculate
\begin{align*}
f_1'( v ) = A v ,\qquad f_2'( v ) =  f .
\end{align*}
Now, using the sum and the chain rules for generalized subgradient (c.f. Propositions~3.35 and 3.37 in $\cite{MOS}$), we obtain
\begin{align*}
 \partial_2\mathcal{L}(w,v) &= f_1'(v)-f_2'(v)+\partial_2(J\circ\gamma)(\gamma w, v)\\
   &\subseteq A  v  -  f  + \gamma^* \partial_2 J(\gamma  w , \gamma  v ),
\end{align*}
which concludes $(ii)$.

\medskip
\noindent
In order to prove $(iii)$, let us fix $ w ,  v _i \in V$ with $i=1, 2$. We take $ \zeta  _i \in \partial_2 \mathcal{L}( w ,  v _i)$. From $(ii)$ there exist $ z _i \in~\partial_2J(\gamma  w , \gamma  v _i)$ such that
\begin{align*}
\zeta_i = A v_i - f + \gamma^* z_i.
\end{align*}
Hence, using $H(A)$(c) and (\ref{RM}), we obtain
\begin{align*}
  &\langle \zeta  _1 -  \zeta  _2,  v _1 -  v _2 \rangle_{V^*\times V}\nonumber\\
  &\quad =\langle A v _1 - A v _2,   v _1 -   v _2  \rangle_{V^*\times V} + \langle \gamma^*  z _1 - \gamma^*  z _2,  v _1 -  v _2  \rangle_{V^*\times V}\nonumber\\
  &\quad\geq m_A\| v _1 -  v _2\|_V^2 + \langle  z _1 -  z _2, \gamma v _1 - \gamma v _2  \rangle_{X^*\times X} \\
  &\quad\geq m_A\| v _1 -  v _2\|_V^2 - m_\alpha\|\gamma v _1 - \gamma v _2\|_X^2\\
  &\quad\geq (m_A - m_\alpha c_\gamma^2)\| v _1 -  v _2\|_V^2.
\end{align*}
From $(H_s)$ we see that $\partial_2 \mathcal{L}( w , \cdot)$ is strongly monotone for every $ w  \in V$. This is equivalent to the fact that $\mathcal{L}( w , \cdot)$ is strongly convex for every $ w  \in V$ (see Theorem~3.4 in $\cite{FLG}$), which implies that it is strictly convex.
\end{pfs}

\medskip
\noindent
The problem under consideration reads as follows.

\medskip
\noindent
\textbf{Problem $\bm{P_{opt}}$:} {\it \ Find $ u  \in V$ such that}
\begin{align*}
  0 \in \partial_2 \mathcal{L}( u , u ).
\end{align*}
We are now in a position to prove the existence and uniqueness result for the above optimization problem.
\begin{lemma} \label{O}
Assume that $H(A)$, $H(J)$, $H(f)$ and $(H_s)$ hold. Then Problem~$P_{opt}$ has a unique solution $ u  \in V$.
\end{lemma}

\begin{pfs}
We introduce operator $\Lambda \colon V \rightarrow V$ defined for all $ w  \in V$ as follows
\begin{align*}
  \Lambda w  = \argmin_{ v  \in V} \mathcal{L}( w , v ).
\end{align*}
From Lemma $\ref{Lprop}$ $(iii)$ we see that operator $\Lambda$ is well defined. Now we prove that the operator $\Lambda$ is a contraction.
Let $ \widehat{u} _i = \Lambda  u _i$ for $ u _i \in V$ fixed, $i=1, 2$. Because of strict convexity of  $\mathcal{L}( w , \cdot)$ we have
\begin{align*}
  \widehat{u} _i =  \argmin_{ v  \in V} \mathcal{L}( u _i,  v ) \quad\mbox{if and only if}\quad 0 \in \partial_2 \mathcal{L}( u _i,  \widehat{u} _i)
\end{align*}
(see Theorem 1.23 in $\cite{HMP}$).
From similar arguments to those used in proofs of Lemmata~$\ref{I}$ and~$\ref{Lprop}$ with fixed first argument of operator $\mathcal{L}$, we have for all $  v  \in V$
\begin{align*}
  &\langle f- A \widehat{u} _i,  v  \rangle_{V^*\times V} \leq J_2^0(\gamma  u _i, \gamma \widehat{u} _i; \gamma v ).
\end{align*}
Taking for $i=1$ value $ v  =  \widehat{u}_2 -  \widehat{u} _1$, for $i=2$ value $ v  =  \widehat{u}_1 -  \widehat{u}_2$ and adding these inequalities, we obtain
\begin{align*}
  &\langle A \widehat{u} _1 - A \widehat{u} _2,  \widehat{u} _1 -  \widehat{u} _2 \rangle_{V^*\times V} \\
  &\leq J_2^0(\gamma  u _1, \gamma \widehat{u} _1; \gamma \widehat{u} _2 - \gamma \widehat{u} _1) + J_2^0(\gamma  u _2, \gamma \widehat{u} _2; \gamma \widehat{u} _1 - \gamma \widehat{u} _2).
\end{align*}
From assumptions $H(A)$(c) and $H(J)$(c), we get
\begin{align*}
  &m_A\| \widehat{u} _1 -  \widehat{u} _2\|_V^2 \leq m_\alpha\|\gamma \widehat{u} _1 - \gamma  \widehat{u} _2\|_X^2 +  m_L\|\gamma  u _1 - \gamma  u _2\|_X\|\gamma \widehat{u} _1 - \gamma \widehat{u} _2\|_X.
\end{align*}
Using the elementary inequality $ab \leq \frac{a^2}{2} + \frac{b^2}{2}$, we obtain
\begin{align*}
  &m_A\| \widehat{u} _1 -  \widehat{u} _2\|_V^2 \leq m_\alpha c^2_\gamma\| \widehat{u} _1 -  \widehat{u} _2\|_V^2 +  \frac{m_Lc_\gamma^2}{2}(\| u _1 -  u _2\|_V^2 + \| \widehat{u} _1 -  \widehat{u} _2\|_V^2).
\end{align*}
Because of $(H_s)$, we can rearrange these terms to get
\begin{align*}
  &\| \widehat{u} _1 -  \widehat{u} _2\|_V^2 \leq \frac{m_L c_\gamma^2}{2m_A - 2m_\alpha c_\gamma^2 - m_L c_\gamma^2}\| u _1 -  u _2\|_V^2.
\end{align*}
Using assumption $(H_s)$ once more, we obtain that the operator $\Lambda$ is a contraction. From the Banach fixed point theorem we know that there exists a unique $ u ^* \in V$ such that $\Lambda u ^* =  u ^*$, so
 $ 0 \in \partial_2 \mathcal{L}( u ^*,  u ^*)$.
\end{pfs}

\noindent
Let us conclude the results from Lemmata~$\ref{I}$, $\ref{Lprop}$ and $\ref{O}$ in the following theorem.

\begin{theorem} \label{I=O}
Assume that $H(A)$, $H(J)$, $H(f)$ and $(H_s)$ hold. Then Problems~$P_{incl}$ and $P_{opt}$ are equivalent,
they have a unique solution $u\in V$ and this solution satisfies
\begin{align*}
  &\| u \|_V \leq c(1+\| f \|_{V^*})
\end{align*}
with a positive constant $c$.
\end{theorem}

\begin{pfs}
Lemma~$\ref{Lprop}$ $(ii)$ implies that every solution to Problem $P_{opt}$ solves Problem~$P_{incl}$. Using this fact, Lemmata~$\ref{I}$ and $\ref{O}$ we see that a unique solution to Problem~$P_{opt}$ is also a unique solution to Problem~$P_{incl}$. Because of the uniqueness of the solution to Problem~$P_{incl}$ we get that Problems~$P_{incl}$ and $P_{opt}$ are equivalent. The estimation in the statement of the theorem follows from Lemma~$\ref{I}$.
\end{pfs}

\section{Numerical scheme} \label{ns}
\noindent
Let $V^h \subset V$ be a family of finite dimensional subspaces with a discretization parameter $h>0$. We present the following discrete scheme of Problem~$P_{opt}$.

\bigskip
\noindent
\textbf{Problem $\bm{P_{opt}^{h}}$:} {\it \ Find  $u^h \in V^h$ such that}
\begin{align*}
0 \in \partial_2 \mathcal{L}(u^h,  u^h).
\end{align*}
We remark that existence of a unique solution to Problem~$P_{opt}^h$ and equivalence to the discrete version of Problem~$P_{incl}$ follow from application of Theorem~$\ref{I=O}$ in this new setting. Now let us present the following main theorem concerning error estimation of introduced numerical scheme.

\begin{theorem} \label{estimation}
Assume that $H(A)$, $H(J)$, $H(f)$ and $(H_s)$ hold. Then for the unique solutions $u$ and $u^h$ to Problems~$P_{opt}$ and $P_{opt}^h$, respectively,  there exists a constant $c>0$ such that
\begin{equation} \label{thminequality}
  \|  u  -  u ^h\|_V^2 \leq c\,\inf\limits_{ v ^h \in V^h}  \Big\{ \| u  -  v ^h \|_V^2 + \|\gamma  u  - \gamma v ^h \|_X + R( u ,  v ^h) \Big\},
\end{equation}
 where a residual quantity is given by
\begin{equation} \label{R}
  R( u ,  v ^h) =  \langle A  u ,  v ^h -  u  \rangle_{V^*\times V} + \langle  f ,  u  -  v ^h \rangle_{V^*\times V}.
\end{equation}
\end{theorem}

\begin{pfs}
\noindent
Let $u$ be a solution to Problem~$P_{opt}$ and $u^h$ be a solution to Problem~$P_{opt}^h$.
Then they are solutions to corresponding inclusion problems and satisfy respectively
\begin{align}
&\langle f - Au, v  \rangle_{V^*\times V}  \leq J_2^0(\gamma u, \gamma u; \gamma v) \quad \mbox{for \ all\ } v \in V, \label{No1}\\[2mm]
&\langle f - Au^h, v  \rangle_{V^*\times V}  \leq J_2^0(\gamma u^h, \gamma u^h; \gamma v) \quad \mbox{for \ all\ } v \in V^h. \label{No2}
\end{align}
Taking (\ref{No1}) with $v=u^h-u$, and (\ref{No2}) with $v=v^h-u^h$, then adding these inequalities, we obtain for all $ v^h \in V^h$
\begin{align}
  &\langle f, v ^h - u \rangle_{V^*\times V} + \langle Au^h - Au, u^h - u \rangle_{V^*\times V} - \langle A u ^h, v ^h -  u  \rangle_{V^*\times V} \nonumber\\
  &\qquad \leq J_2^0(\gamma u , \gamma u ; \gamma u ^h - \gamma u ) + J_2^0(\gamma  u ^h, \gamma  u ^h; \gamma  v ^h - \gamma  u ^h). \label{27}
\end{align}
We observe that by subadditivity of generalized directional derivative (cf. \cite{MOS}, Proposition 3.23(i)) and $H(J)$(c), we have
\begin{align}
  & J_2^0(\gamma u, \gamma u; \gamma u ^h - \gamma u) + J_2^0(\gamma u^h, \gamma u^h; \gamma v^h - \gamma u^h )\nonumber\\
  &\leq J_2^0(\gamma u, \gamma u; \gamma u^h - \gamma u) + J_2^0(\gamma u^h, \gamma u^h; \gamma u - \gamma u^h) + J_2^0(\gamma u^h, \gamma u^h; \gamma v^h - \gamma u)\nonumber\\
  &\leq (m_\alpha + m_L) \|\gamma u^h -\gamma u \|_X^2 + \left( c_0 + (c_1+c_2) \|\gamma u^h\|_X \right)\| \gamma  v ^h - \gamma  u \|_X. \label{4}
\end{align}
From the statement of Lemma $\ref{I}$ applied to discrete version of Problem $P_{incl}$ we get that $\|\gamma u^h\|_X \leq c_{\gamma}\|u^h\|_V\leq c\,(1+\|f\|_{V^*})$ is uniformly bounded with respect to $h$. Hence, returning to (\ref{27}) and using ($\ref{4}$), we obtain for all $ v ^h \in V^h$
 \begin{align*}
  &\langle Au^h - Au, u^h - u \rangle_{V^*\times V} \leq  \langle A u ^h - Au, v ^h -  u  \rangle_{V^*\times V} + \langle A u, v ^h -  u  \rangle_{V^*\times V} \nonumber\\
  &\qquad + \langle f, u-v^h \rangle_{V^*\times V} + (m_\alpha + m_L)c^2_\gamma \|u^h - u \|_V^2 + c \, \|\gamma v^h  - \gamma u\|_X .
\end{align*}
By assumption $H(A)$ and definition~(\ref{R}), we get for all $v^h\in V^h$
\begin{align*}
  &m_A \|  u ^h -  u \|_V^2 \leq c \, \| u ^h -  u \|_V\| v ^h -  u \|_V + R( u ,  v ^h) \nonumber\\
   &\qquad +  (m_\alpha + m_L) c_\gamma^2\|  u  -  u ^h \|_V^2 + c \, \| \gamma  u  - \gamma  v ^h \|_X. 
\end{align*}
Finally, the elementary inequality $ab\leq \varepsilon a^2 + \frac{b^2}{4\varepsilon}$ with $\varepsilon > 0$ yields
\begin{align*}
  & m_A \| u  -  u ^h\|_V^2 \leq \varepsilon\|  u  -  u ^h\|_V^2 + \frac{c^2}{4\varepsilon}\|  u  -  v ^h\|_V^2 + R( u ,  v ^h) \\
  &\qquad + (m_\alpha + m_L) c_\gamma^2\|  u  -  u ^h \|_V^2 + c \, \| \gamma  u  - \gamma  v ^h \|_X.
\end{align*}
This is equivalent for all $ v ^h \in V^h$ to
\begin{align*}
  & \Big(m_A - (m_\alpha + m_L) c_\gamma^2 - \varepsilon\Big) \|  u  -  u ^h\|_V^2 \leq \frac{c}{\varepsilon}\| u  -  v ^h \|_V^2 + R( u ,  v ^h) + c \, \| \gamma  u  - \gamma  v ^h \|_X.
\end{align*}
Taking sufficiently small $\varepsilon$ and using $(H_s)$ we obtain the desired conclusion.
\end{pfs}

\section{Application to Contact Mechanics} \label{mcp}
\noindent
In this section we apply the results of previous sections to a sample mechanical contact problem.
Let us start by introducing the physical setting and notation useful in the problem.

An elastic body occupies a domain $\Omega \subset \mathbb{R}^{d}$, where $d = 2, 3$ in application. We assume that its boundary $\Gamma$ is divided into three disjoint measurable parts:
$\Gamma_{D}, \Gamma_{C}, \Gamma_{N}$, where the part $\Gamma_D$ has a positive measure.
Additionally $\Gamma$ is Lipschitz continuous, and therefore the outside normal vector $\bm{\nu}$ to $\Gamma$ exists a.e. on the boundary.
The body is clamped on $\Gamma_{D}$, i.e. its displacement is equal to $\bm{0}$ on this part of boundary.
A surface force of density $\bm{f}_N$ acts on the boundary~$\Gamma_{N}$ and a body force of density $\bm{f}_0$ acts in $\Omega$. The contact phenomenon on $\Gamma_{C}$ is modeled using general subdifferential inclusions. We are interested in finding the displacement of the body in a static state.

Let us denote by ``$\cdot$'' and $\|\cdot\|$ the scalar product and the Euclidean norm in $\mathbb{R}^{d}$ or $\mathbb{S}^{d}$, respectively, where  $\mathbb{S}^{d} = \mathbb{R}^{d \times d}_{sym}$. Indices $i$ and $j$ run from $1$ to $d$ and the index after a comma represents the partial derivative with respect to the corresponding component of the independent variable.
Summation over repeated indices is implied. We denote the divergence operator by $\textrm{Div }\bm{\sigma} = (\sigma_{ij,j})$. The standard Lebesgue and Sobolev spaces $L^2(\Omega)^d = L^2(\Omega;\mathbb{R}^d)$ and $H^1(\Omega)^d = H^1(\Omega;\mathbb{R}^d)$ are used.
The linearized (small) strain tensor for displacement $\bm{u} \in H^1(\Omega)^d$ is defined by
\begin{equation}\nonumber
 \bm{\varepsilon}(\bm{u})=(\varepsilon_{ij}(\bm{u})), \quad \varepsilon_{ij}(\bm{u}) = \frac{1}{2}(u_{i,j} + u_{j,i}).
\end{equation}
Let $u_\nu= \bm{u}\cdot \bm{\nu}$ and $\sigma_\nu= \bm{\sigma}\bm{\nu} \cdot \bm{\nu}$ be the normal components of $\bm{u}$ and $\bm{\sigma}$, respectively,  and let $\bm{u}_\tau =\bm{u}-u_\nu\bm{\nu}$ and $\bm{\sigma}_\tau =\bm{\sigma}\bm{\nu}-\sigma_\nu\bm{\nu}$ be their tangential components, respectively. In what follows, for simplicity, we sometimes do not indicate explicitly the dependence of various functions on the spatial variable $\bm{x}$.

Now let us introduce the classical formulation of considered mechanical contact problem.

\medskip
\noindent
\textbf{Problem $\bm{P}$:} \textit{Find a displacement field $\bm{u}\colon \Omega \rightarrow \mathbb{R}^{d}$ and a stress field $\bm{\sigma}\colon \Omega \rightarrow \mathbb{S}^{d}$ such that}

\begin{align}
  \bm{\sigma}  = \mathcal{A}(\bm{\varepsilon}(\bm{u})) \qquad &\textrm{ in } \Omega \label{P1}\\
  \textrm{Div }\bm{\sigma} + \bm{f}_{0} = \bm{0}  \qquad &\textrm{ in } \Omega \label{P2}\\
  \bm{u} = \bm{0}  \qquad &\textrm{ on } \Gamma_{D} \label{P3}\\
  \bm{\sigma}\bm{\nu} = \bm{f}_{N} \qquad &\textrm{ on } \Gamma_{N} \label{P4}\\
   -\sigma_{\nu} \in \partial j_{\nu}(u_{\nu}) \qquad &\textrm{ on } \Gamma_{C} \label{P5}\\
   -\bm{\sigma_{\tau}} \in h_{\tau}(u_{\nu})\,\partial j_{\tau}(\bm{u_{\tau}}) \qquad &\textrm{ on } \Gamma_{C} \label{P6}
\end{align}

\noindent
Here, equation~\eqref{P1} represents an elastic constitutive law and $\mathcal{A}$ is an elasticity operator.
Equilibrium equation~(\ref{P2}) reflects the fact that problem is static.
Equation~(\ref{P3}) represents clamped boundary condition on $\Gamma_{D}$ and~(\ref{P4}) represents the action of the traction on $\Gamma_{N}$.
Inclusion~(\ref{P5}) describes the response of the foundation in normal direction, whereas the friction is modeled by inclusion~(\ref{P6}), where $j_\nu$ and $j_\tau$ are given superpotentials, and $h_\tau$ is a~given friction bound.


We consider the following Hilbert spaces
\begin{align*}
 &\mathcal{H} = L^2(\Omega;\mathbb{S}^{d}), \qquad
 V = \{\bm{v} \in H^1(\Omega)^d\ |\ \bm{v} = \bm{0} \textrm{ on } \Gamma_{D}\},
\end{align*}
endowed with the inner scalar products
\begin{align*}
 &(\bm{\sigma},\bm{\tau})_{\mathcal{H}} = \int_\Omega \sigma_{ij}\tau_{ij} \, dx, \qquad
 (\bm{u},\bm{v})_V = (\bm{\varepsilon}(\bm{u}),\bm{\varepsilon}(\bm{v}))_{\mathcal{H}},
\end{align*}
respectively.
The fact that space $V$ equipped with the norm $\|\cdot\|_V$ is complete follows from Korn's inequality, and its application is allowed because we assume that $meas(\Gamma_{D}) > 0$.
We consider the trace operator $\gamma \colon V \to L^2(\Gamma_{C})^d=X$.
By the Sobolev trace theorem we know that $\gamma \in  \mathcal{L}(V, X)$ with the norm equal to $c_\gamma$.

\medskip
\noindent
Now we present the hypotheses on data of Problem~$P$.

\medskip
\noindent
$\underline{H({\mathcal{A}})}:$ \quad ${\mathcal{A}} \colon \Omega \times {\mathbb S}^d \to {\mathbb S}^d$ satisfies
\begin{enumerate}
  \item[(a)]
     $\mathcal{A}(\bm{x},\bm{\tau}) = (a_{ijkh}(\bm{x})\tau_{kh})$
     for all $\bm{\tau} \in {\mathbb S}^d$, a.e. $\bm{x}\in\Omega,\ a_{ijkh} \in L^{\infty}(\Omega),$
  \item[(b)]
    $\mathcal{A}(\bm{x},\bm{\tau}_1) \cdot \bm{\tau}_2 = \bm{\tau}_1 \cdot  \mathcal{A}(\bm{x},\bm{\tau}_2)$ for all $\bm{\tau}_1, \bm{\tau}_2 \in {\mathbb S}^d$, a.e. $\bm{x}\in\Omega$,
  \item[(c)]
     there exists $m_{\mathcal{A}}>0$ such that $\mathcal{A}(\bm{x},\bm{\tau}) \cdot \bm{\tau} \geq m_{\mathcal{A}} \|\bm{\tau}\|^2$ for all $\bm{\tau} \in {\mathbb S}^d$, a.e. $\bm{x}\in\Omega$.
\end{enumerate}

\medskip
\noindent
$\underline{H(j_{\nu})}:$ \quad $j_{\nu} \colon \Gamma_C \times \mathbb{R} \to \mathbb{R}$ satisfies
\begin{enumerate}
  \item[(a)]
    $j_{\nu}(\cdot, \xi)$ is measurable on $\Gamma_C$ for all $\xi \in \mathbb{R}$ and there exists $e \in L^2(\Gamma_C)$ such that $j_{\nu}(\cdot,e(\cdot))\in L^1(\Gamma_C)$,
  \item[(b)]
    $j_{\nu}(\bm{x}, \cdot)$ is locally Lipschitz continuous on $\mathbb{R}$ for a.e. $\bm{x} \in \Gamma_C$,
  \item[(c)]
    there exist $c_{\nu0}, c_{\nu1} \geq 0$ such that \\[2mm]
    \hspace*{1cm}$|\partial_2 j_{\nu}(\bm{x}, \xi)| \leq c_{\nu0} + c_{\nu1}|\xi|$\quad for all $\xi \in \mathbb{R}$, a.e. $\bm{x} \in \Gamma_C$,
  \item[(d)]
    there exists $\alpha_{\nu} \geq 0$ such that \\[2mm]
    \hspace*{1cm}$(j_{\nu})_2^0(\bm{x},\xi_1;\xi_2-\xi_1) + (j_{\nu})_2^0(\bm{x},\xi_2;\xi_1-\xi_2)\leq \alpha_{\nu}|\xi_1-\xi_2|^2$ \\[2mm]
     for all $\xi_1, \xi_2 \in \mathbb{R}$, a.e. $\bm{x} \in \Gamma_C$.
\end{enumerate}

\medskip
\noindent
$\underline{H(j_{\tau})}:$ \quad $j_{\tau} \colon \Gamma_C \times \mathbb{R}^{d} \to \mathbb{R}$ satisfies
\begin{enumerate}
  \item[(a)]
    $j_{\tau}(\cdot, \bm{\xi})$ is measurable on $\Gamma_C$ for all $\bm{\xi} \in \mathbb{R}^{d}$ and there exists $\bm{e} \in L^2(\Gamma_C)^{d}$ such that $j_{\tau}(\cdot,\bm{e}(\cdot))\in L^1(\Gamma_C)$,
  \item[(b)]
    there exists $c_{\tau}>0$ such that \\[2mm]
    \hspace*{1cm}$|j_{\tau}(\bm{x}, \bm{\xi}_1) - j_{\tau}(\bm{x}, \bm{\xi}_2)| \leq c_{\tau} \|\bm{\xi}_1 - \bm{\xi}_2\|$\quad for all $\bm{\xi}_1, \bm{\xi}_2 \in \mathbb{R}^d$, a.e. $\bm{x} \in \Gamma_C$,
  \item[(c)]
    there exists $\alpha_{\tau} \geq 0$ such that  \\[2mm]
   \hspace*{1cm} $(j_{\tau})_2^0(\bm{x},\bm{\xi}_1;\bm{\xi}_2-\bm{\xi}_1) + (j_{\tau})_2^0(\bm{x},\bm{\xi}_2;\bm{\xi}_1-\bm{\xi}_2)\leq \alpha_\tau\|\bm{\xi}_1-\bm{\xi}_2\|^2$\\[2mm]
     for all $\bm{\xi}_1, \bm{\xi}_2 \in \mathbb{R}^{d}$, a.e. $\bm{x} \in \Gamma_C$.
\end{enumerate}

\medskip
\noindent
$\underline{H(h)}:$ \quad $h_{\tau} \colon \Gamma_C \times \mathbb{R} \to \mathbb{R}$ satisfies
\begin{enumerate}
  \item[(a)]
    $h_{\tau}(\cdot, \eta)$ is measurable on $\Gamma_C$ for all $\eta \in \mathbb{R}$,
  \item[(b)]
    there exists $\overline{h_{\tau}} > 0$ such that $0 \leq h_{\tau}(\bm{x}, \eta) \leq \overline{h_{\tau}}$ for all $\eta \in \mathbb{R}$, a.e. $\bm{x} \in \Gamma_C$,
  \item[(c)]
    there exists $L_{h_\tau}>0$ such that \\[2mm]
    \hspace*{1cm} $|h_{\tau}(\bm{x}, \eta_1) - h_{\tau}(\bm{x}, \eta_2)| \leq L_{h_\tau} |\eta_1 - \eta_2|$\quad for all $\eta_1, \eta_2 \in \mathbb{R}$, a.e. $\bm{x} \in \Gamma_C$.
\end{enumerate}

\medskip
\noindent
$(\underline{H_0}): \quad \bm{f}_0 \in L^2(\Omega)^d, \quad \bm{f}_N \in L^2(\Gamma_N)^d$.

\bigskip
We remark that condition $H(j_{\tau})$(b) is equivalent to the fact that $j_{\tau}(\bm{x},\cdot)$ is locally Lipschitz continuous and there exists $c_{\tau} \geq 0$ such that $\|\partial_2 j_{\tau}(\bm{x} ,  \bm{\xi} )\| \leq c_{\tau}$ for all $\bm{\xi} \in \mathbb{R}^{d}$ and a.e. $\bm{x}  \in \Gamma_C$.

\medskip
\noindent
Using the standard procedure, the Green formula and the definition of generalized subdifferential, we obtain a weak formulation of Problem~$P$ in the form of hemivariational inequality.

\medskip
\noindent
\textbf{Problem $\bm{P_{hvi}}$:} {\it \ Find a displacement $\bm{u} \in V$ such that for all $\bm{v} \in V$}
\begin{align}
  &\langle A\bm{u}, \bm{v} \rangle_{V^*\times V}  +\int_{\Gamma_C} j_3^0 (\bm{x}, \gamma \bm{u}(\bm{x}), \gamma \bm{u}(\bm{x}); \gamma \bm{v}(\bm{x}))\, da \geq \langle \bm{f}, \bm{v} \rangle_{V^*\times V}. \label{PV}
\end{align}

\medskip
\noindent
Here, the operator $A \colon V \to V^*$ and $\bm{f} \in V^*$ are defined for all $\bm{w},\bm{v} \in V$ as follows
\begin{align*}
 &\langle A\bm{w}, \bm{v} \rangle_{V^*\times V} = (\mathcal{A}(\bm{\varepsilon}(\bm{w})),\bm{\varepsilon}(\bm{v}))_{\mathcal{H}}, \nonumber\\
 &\langle \bm{f}, \bm{v} \rangle_{V^* \times V} = \int_{\Omega}\bm{f}_{0}\cdot \bm{v}\, dx + \int_{\Gamma_{N}}\bm{f}_{N}\cdot \gamma \bm{v}\, da\nonumber
\end{align*}
and $j\colon \Gamma_C\times \mathbb{R}^d \times \mathbb{R}^d \to \mathbb{R}$ is defined for all $\bm{\eta}, \bm{\xi} \in \mathbb{R}^d$ and $\bm{x}\in \Gamma_C$ by
\begin{align}
j(\bm{x}, \bm{\eta}, \bm{\xi}) =  j_{\nu}(\bm{x}, \xi_{\nu}) + h_{\tau}(\bm{x}, \eta_{\nu})\, j_{\tau}(\bm{x}, \bm{\xi}_{\tau}). \label{Jj}
\end{align}
It is easy to check that under assumptions $H(\mathcal{A})$ and $(H_0)$, the operator~$A$ and the functional~$\bm{f}$ satisfy $H(A)$ and $H(f)$, respectively. 
We also define the functional $J \colon L^2(\Gamma_C)^d \times L^2(\Gamma_C)^d \to \mathbb{R}$ for all $\bm{w}, \bm{v} \in L^2(\Gamma_C)^d $ by
\begin{align}
J(\bm{w}, \bm{v}) =  \int_{\Gamma_C} j(\bm{x}, \bm{w}(\bm{x}), \bm{v}(\bm{x}))\, da,\label{J}
\end{align}

\medskip
\noindent
Below we present some properties of the functional~$J$. The proof of the following lemma is similar to the proof of Corollary~4.15 in $\cite{MOS}$ and is skipped here.

\begin{lemma} \label{LJ}
Assumptions $H(j_{\nu})$, $H(j_{\tau})$ and $H(h)$ imply that functional $J$ defined by {\rm{(\ref{Jj})-(\ref{J})}} satisfies $H(J)$.
\end{lemma}

\medskip
\noindent
With the above properties, we have the following existence and uniqueness result for Problem~$P_{hvi}$.
\begin{theorem} \label{O=H}
If assumptions $H(\mathcal{A})$, $H(j_{\nu})$, $H(j_{\tau})$, $H(h)$, $(H_0)$ and $(H_s)$ hold, then Problems~$P_{hvi}$ and $P_{incl}$ are equivalent. Moreover, they have a unique solution~$\bm{u}\in V$ and this solution satisfies
\begin{align*}
  &\| \bm{u} \|_V \leq c\, (1+\| \bm{f} \|_{V^*})
\end{align*}
with a positive constant $c$.
\end{theorem}

\begin{pfs}
We notice that the assumptions of Theorem~\ref{I=O} are satisfied. This implies that Problem~$P_{incl}$ has a unique solution. By (\ref{9}) and Corollary~4.15~(iii) in $\cite{MOS}$ we get that every solution to Problem $P_{incl}$ solves Problem $P_{hvi}$. Using similar technique as in the proof of Lemma $\ref{I}$, we can show that if Problem $P_{hvi}$  has a solution, it is unique. Combining these facts we obtain our assertion.
\end{pfs}

\noindent
We conclude this section by providing a sample error estimate under additional assumptions on the solution regularity. We consider a polygonal domain $\Omega$ and a space of continuous piecewise affine functions $V^h$. We introduce the following discretized version of Problem $P_{hvi}$.

\medskip
\noindent
\textbf{Problem $\bm{P_{hvi}^h}$:} \textit{Find a displacement $\bm{u}^h \in V^h$ such that for all $\bm{v}^h \in V^h$}
\begin{align}
  &\langle A\bm{u}^h, \bm{v}^h \rangle_{V^*\times V}  + \int_{\Gamma_{C}} j_3^0(\bm{x},\gamma \bm{u}^h(\bm{x}), \gamma \bm{u}^h(\bm{x}); \gamma \bm{v}^h(\bm{x}))\, da \geq \langle \bm{f}, \bm{v}^h \rangle_{V^*\times V}.
\end{align}

\begin{theorem}
Assume $H(\mathcal{A})$, $H(j_{\nu})$, $H(j_{\tau})$, $H(h)$, $(H_0)$ and $(H_s)$ and assume the solution regularity $\bm{u} \in H^2(\Omega)^d$,  $\bm{\sigma}\bm{\nu} \in L^2(\Gamma_C)^d$. Then, for the solution $\bm{u}$ to Problem~$P_{hvi}$ and the solution $\bm{u}^h$ to Problem~$P_{hvi}^h$ there exists a constant $c>0$ such that 
\begin{equation} \nonumber
  \|\bm{u} - \bm{u}^{h} \|_V \leq c\, h.
\end{equation}
\end{theorem}

\begin{pfs}
We denote by $\Pi^h \bm{u} \in V^h$ the finite element interpolant of $u$. By the standard finite element interpolation error bounds (see \cite{Cia}) we have for all $\bm{\eta} \in H^2(\Omega)^d$
\begin{align}
  &\|\bm{\eta} - \Pi^h\bm{\eta}\|_V \leq c\,h\,\|\bm{\eta}\|_{H^2(\Omega)^d},\label{ie1}\\
  &\|\gamma \bm{\eta} - \gamma \Pi^h\bm{\eta}\|_{L^2(\Gamma_C)^d} \leq c\,h^2\,\|\bm{\eta}\|_{H^2(\Gamma_C)^d}\label{ie2}.
\end{align} 
We now bound the residual term defined by (\ref{R}) using similar procedure to one described in \cite{H}. Let $\bm{v} = \pm \bm{w}$ in inequality (\ref{PV}), where the arbitrary function $\bm{w}\in V$ is such that $\bm{w} \in C^\infty(\overline\Omega)^d$ and $\bm{w} = \bm{0}$ on $\Gamma_D \cup \Gamma_C$. Then we obtain the identity
\begin{align*}
  &\langle A\bm{u}, \bm{w} \rangle_{V^*\times V}  = \langle \bm{f}, \bm{w} \rangle_{V^*\times V}.
\end{align*}
From this identity, using fundamental lemma of calculus of variations, we can deduce that
\begin{align}
    \textrm{Div }\mathcal{A}(\bm{\varepsilon}(\bm{u})) + \bm{f}_{0} = \bm{0}  \qquad &\textrm{ in } \Omega, \label{d1}\\
    \bm{\sigma}\bm{\nu} = \bm{f}_{N} \qquad &\textrm{ on } \Gamma_{N}. \label{d2}
\end{align}  
We multiply equation (\ref{d1}) by $\bm{v}^h - \bm{u}$ and obtain
\begin{align}
    \int_{\Gamma}  \bm{\sigma}\bm{\nu} \cdot (\gamma\bm{v}^h - \gamma\bm{u})\, da - \int_{\Omega} \mathcal{A}(\bm{\varepsilon}(\bm{u}))  \cdot \bm{\varepsilon}(\bm{v}^h - \bm{u})\, dx + \int_{\Omega} \bm{f}_{0} \cdot (\bm{v}^h - \bm{u})\, dx = 0.
\end{align}
Using the homogenous Dirichlet boundary condition of $\bm{v}^h - \bm{u}$ on $\Gamma_D$ and the traction boundary condition given by (\ref{d2}) we have
\begin{align}
    \langle A \bm{u} ,  \bm{v}^h - \bm{u}  \rangle_{V^*\times V} =  \int_{\Gamma_{C}} \bm{\sigma}\bm{\nu}\cdot (\gamma\bm{v}^h - \gamma\bm{u})\, da  + \langle \bm{f} ,  \bm{v}^h - \bm{u} \rangle_{V^*\times V}.
\end{align}
Using this and ($\ref{R}$) we obtain
\begin{align} 
 &R(\bm{u}, \bm{v}^{h}) =  \int_{\Gamma_{C}} \bm{\sigma}\bm{\nu}\cdot (\gamma\bm{v}^h - \gamma\bm{u})\, da \leq c\,\|\gamma\bm{u} - \gamma\bm{v}^h\|_{L^2(\Gamma_C)^d}. \label{R2}
\end{align}
From inequalities (\ref{thminequality}), (\ref{ie1}), (\ref{ie2}) and (\ref{R2})  we get
\begin{align*} \nonumber
  &\|\bm{u} - \bm{u}^{h} \|_V^{2} \leq c \, \Big(\| \bm{u}  -  \Pi^h\bm{u} \|_V^2 + \| \gamma \bm{u}  - \gamma \Pi^h\bm{u} \|_{L^2(\Gamma_C)^d} \Big)
    \leq c\,h^2,
\end{align*}
and we obtain required estimation.
\end{pfs}

\section{Simulations} \label{s}
\noindent In this section we present results of our computational simulations. From Theorems~\ref{I=O} and $\ref{O=H}$ we know that Problems~$P_{hvi}$ and $P_{opt}$ are equivalent. Hence, we can apply numerical scheme~$P_{opt}^h$ and use Theorem $\ref{estimation}$ to approximate solution of $P_{hvi}$.
We employ Finite Element Method and use space $V^h$ of continuous piecewise affine functions as a family of approximating subspaces.
The idea for algorithm used to calculate solution of discretized problem is based on the proof of Lemma $\ref{O}$ and is described by Algorithm~1.

\begin{algorithm}[ht]
\caption{Iterative optimization algorithm}
\begin{algorithmic}
\medskip
\State Let $\varepsilon > 0$ and $\bm{u}^h_0$ be given
\State $k\gets 0$
\Repeat
    \State $k\gets k+1$
    \State $\bm{u}^h_{k} =  \argmin_{\bm{v}^h \in V^h} \mathcal{L}(\bm{u}^h_{k-1}, \bm{v}^h)$
\Until {$\|\bm{u}^h_{k} - \bm{u}^h_{k-1}\|_V \leq \varepsilon$}\\
\Return $\bm{u}^h_{k}$
\medskip
\end{algorithmic}
\end{algorithm}

In order to minimize not necessarily differentiable function $\mathcal{L}(\bm{w}^h, \cdot)$ we use Powell's conjugate direction method. For a starting point $\bm{u}_0^h$ we take a solution to problem with $\bm{\sigma}\bm{\nu}=\bm{0}$ on $\Gamma_{C}$, although it can be chosen arbitrarily.

\medskip

We set $d=2$ and consider a rectangular set $\Omega = [0,2] \times [0,1]$ with following parts of the boundary
\begin{align*}
 &\Gamma_{D} = \{0\} \times [0,1], \quad \Gamma_{N} = ([0,2] \times \{1\}) \cup (\{2\} \times [0,1]), \quad \Gamma_{C} = [0,2] \times \{0\}.
\end{align*}

\noindent The elasticity operator $\mathcal{A}$ is defined by
\begin{align*}
 &\mathcal{A}(\bm{\tau}) = 2\eta\bm{\tau} + \lambda \mbox{tr}(\bm{\tau})I,\qquad \bm{\tau} \in \mathbb{S}^2.
\end{align*}

\noindent Here $I$ denotes the identity matrix, $\mbox{tr}$ denotes the trace of the matrix, $\lambda$ and $\eta$ are the Lame coefficients, $\lambda, \eta >0$.
In our simulations we take the following data 

\begin{align*}
 &\lambda = \eta = 4, \\
 &\bm{u}_{0}(\bm{x}) = (0,0), \quad \bm{x} \in \Omega,\\
 &j_\nu(\bm{x}, \xi)= \left \{ \begin{array}{ll}
   0, \quad \xi \in (-\infty,\, 0), \\
   10\, \xi^2, \quad \xi \in [0,\, 0.1), \\
   0.1, \quad \xi \in [0.1,\, \infty), \\
  \end{array} \right. \bm{x} \in \Gamma_C,\\
 &j_{\tau}(\bm{x}, \bm{\xi}) = \ln (\|\bm{\xi}\| + 1), \quad \bm{\xi} \in \mathbb{R}^2,\ \bm{x} \in \Gamma_C,\\
 &h_{\tau}(\bm{x}, \eta) = \left \{ \begin{array}{ll}
   0, \quad \eta \in (-\infty, 0), \\
   8\, \eta, \quad \eta \in [0,\infty), \\
  \end{array} \right. \bm{x} \in \Gamma_C,\\
 &\bm{f}_0(\bm{x}) = (-1.2,\, -0.9), \quad \bm{x} \in \Omega,\\
 &\bm{f}_N(\bm{x}) = (0,0), \quad \bm{x} \in \Omega.
\end{align*}

\begin{figure}[ht]
\centering
    \includegraphics[width=0.6\linewidth]{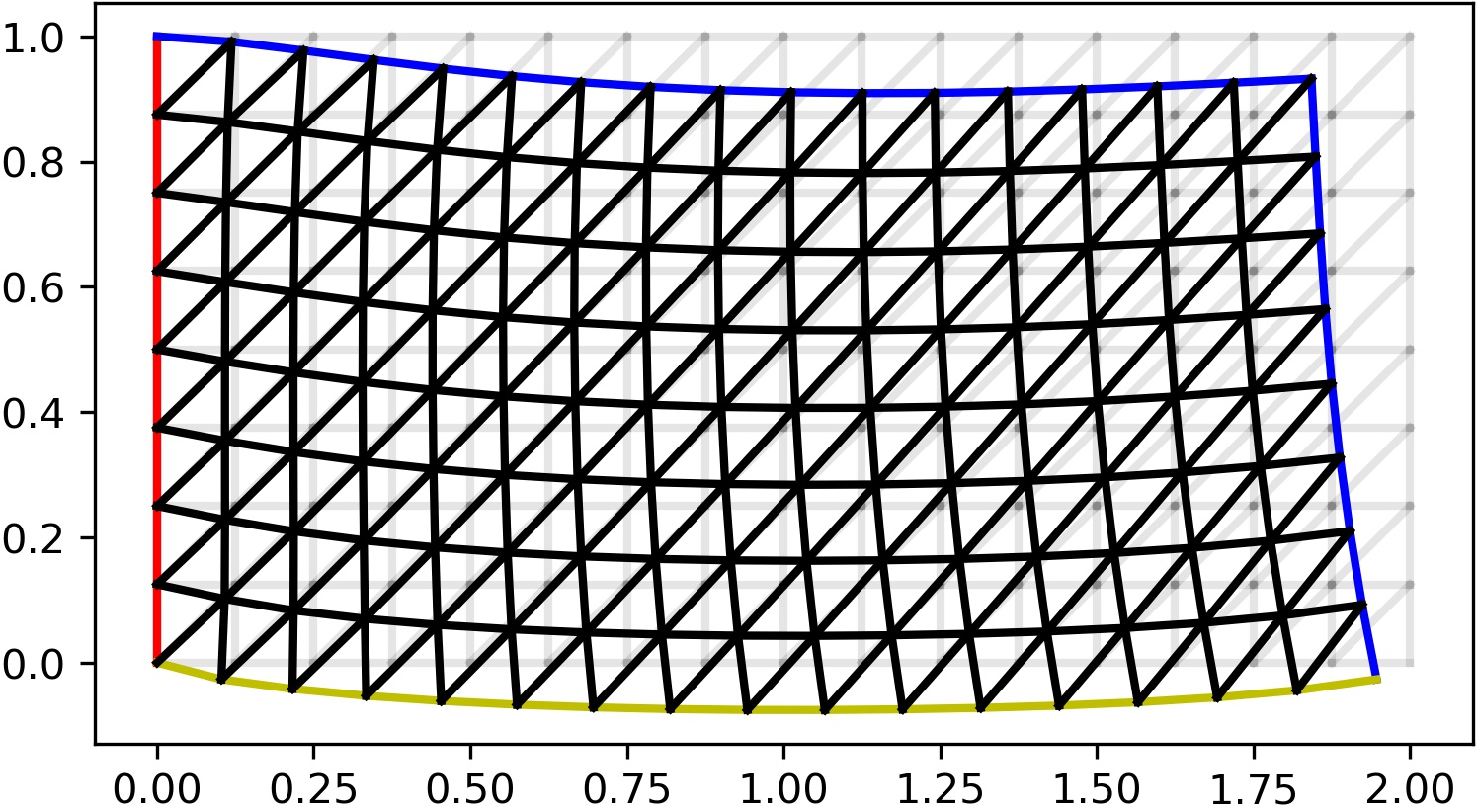}
    \caption{The behaviour of the body} \label{figOne}
\end{figure}

\noindent Both functions $j_{\nu}$ and $j_{\tau}$ are nondifferentiable and nonconvex. In Figure \ref{figOne} we present output obtained for chosen data. We push the body down and to the left with a force $\bm{f}_0$. As a result the body penetrates the foundation, but because of frictional forces it is squeezed to the left more in the higher part than in the lower part. 

\begin{figure}[H]
\centering
  \includegraphics[width=0.6\linewidth]{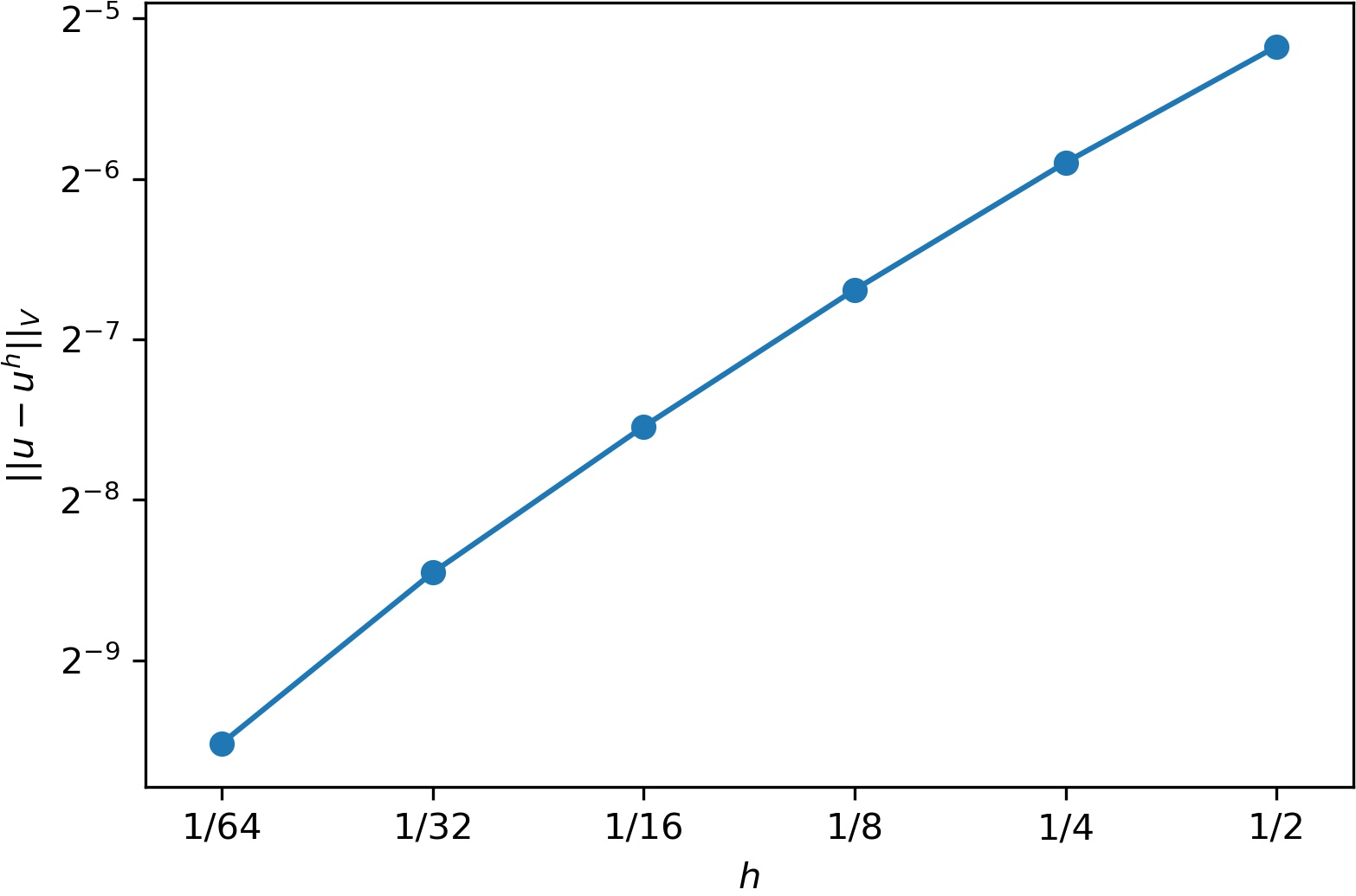}
  \caption{Numerical errors} \label{figFive}
\end{figure}

In order to illustrate the error estimate obtained in Section \ref{mcp}, we present a comparison of numerical errors $\|\bm{u} - \bm{u}^h\|_V$  computed for a sequence of solutions to discretized problems. We use a uniform discretization of the problem domain according to the spatial discretization parameter $h$. The boundary $\Gamma_C$ of $\Omega$ is divided into $1/h$ equal parts. We start with $h = 1$, which is successively halved. The numerical solution corresponding to $h = 1/128$ was taken as the “exact” solution $\bm{u}$. The numerical results
are presented in Figure \ref{figFive}, where the dependence of the error estimate $\|\bm{u}  - \bm{u}^h\|_V$ with respect to $h$ is plotted on a log-log scale. A first order convergence can be observed, providing numerical evidence of the theoretical optimal order error estimate obtained at the end of Section \ref{mcp}.

\medskip
\noindent {\bf Acknowledgments}\\
The project leading to this application has received funding from the European Union's Horizon 2020 research and innovation programme under the Marie Sklo\-do\-wska-Curie grant agreement no. 823731 CONMECH,
National Science Center of Poland under Maestro Project no. UMO-2012/06/A/ST1/00262.



\end{document}